\documentclass[12pt]{amsart}
\usepackage{amsmath}
\usepackage{amssymb}

\newtheorem{lemma}{Lemma}[section]
\newtheorem{theorem}[lemma]{Theorem}
\newtheorem{corollary}[lemma]{Corollary}

\newtheorem{example}[lemma]{Example}
\newtheorem{proposition}[lemma]{Proposition}

\def\pf{\noindent{\bf Proof.  }}
\def\rank{\mathrm{rank}\,}

\begin{document}
\date{}
\title{The expression of  Moore--Penrose inverse of $A-XY^*$}
\author[Fapeng Du, Yifeng Xue]{Fapeng Du, Yifeng Xue $^*$}
\address{Fapeng Du, email: jsdfp@163.com\newline
Department of mathematics, East China Normal University\newline
  Shanghai 200241, P.R. China and \newline
  Collage of mathematics and Physical Sciences, Xuzhou Institute of Technology\newline
  Xuzhou 221008, Jiangsu Province, P.R. China\newline
  $^*$ Department of mathematics, East China Normal University\newline
  Shanghai 200241, P.R. China}
  \thanks{$^*$ Corresponding author; email: yfxue@math.ecnu.edu.cn}
\thanks{Project supported by
Natural Science Foundation of China (no.10771069) and Shanghai Leading Academic Discipline
Project(no.B407)}

\maketitle
\baselineskip 18pt
\begin{abstract}
Let $K,\,H$ be Hilbert spaces and let $L(K,H)$ denote the set of all bounded linear operators
from $K$ to $H$.
Let $A \in L(H)\triangleq L(H,H)$ with $R(A)$ closed and $X,Y \in L(K,H)$ with $R(X)\subseteq
R(A),R(Y)\subseteq R(A^*)$. In this short note, we give some new expressions of the Moore--Penrose
inverse  $(A-XY^*)^+$ of $A-XY^*$ under certain suitable conditions.
\end{abstract}
\vspace{1mm}

\section{Introduction}

Let $A$ be a nonsingular $m \times m$ matrix and $X,\,Y$ be two $m \times n$ matrices.  It is known that $A-XY^*$
 is nonsingular iff $I_n - Y^*A^{-1}X$ is nonsingular, and in which case the well known
 Shermen-Morrison-Woodbury formula (SMW) can be expressed as
\begin{equation}\label{EQ}
(A-XY^*)^{-1}=A^{-1}+A^{-1}X(I_n-Y^*A^{-1}X)^{-1}Y^*A^{-1}
\end {equation}
This formula and some related formula have a lot of applications in statistics, networks,
optimization and partial differential equations. Please see \cite{HS,KR,WWH} for details.
Clearly, the formula  (\ref{EQ}) fails when $A$ or $A-XY^*$ is singular. Steerneman and
Kleij in \cite{TFP}  proved that when A is singular and
$I_n - Y^*A^+X$ is nonsingular, then
\begin{equation*}
(A-XY^*)^{+}=A^{+}+A^{+}X(I_n-Y^*A^+X)^{-1}Y^*A^{+}
\end{equation*}
under conditions that
$$
\rank(A,X)=\rank A,\quad\rank\begin{pmatrix}A\\ Y^*\end{pmatrix}=\rank A.
$$
He also showed that if $A$ is nonsingular and $Y^*A^{-1}X=I_n$, then
\begin{equation}\label{EQ2}
(A-XY^*)^{+}=(I_m -X_1X_1^+)A^{-1}(I_m-Y_1Y_1^+)
\end{equation}\label{FQ}
where $X_1=A^{-1}X ,~~ Y_1=(A^{-1})^*Y$ (cf. \cite[Theorem 3]{TFP}).

Recently Chen, Hu and Xu studied the Moore-Penrose inverse of $A-XY^*$ when $A \in L(H)$ and $X,Y \in L(K,H)$ in \cite{CHX}.
They prove that if $A$ is invertible and $A-XY^*,~X,~Y$ have closed ranges, then
\begin{equation*}
(A-XY^*)^{+}=(I -X_1X_1^+)A^{-1}(I-Y_1Y_1^+)
\end{equation*}
iff $Y_1^*XY_1^*=Y_1^*,\ XY_1^*X=X$,where $X_1=A^{-1}X,~~ Y_1=(A^{-1})^*Y$.
This result generalizes Theorem 3 of \cite{TFP}.

In this paper we assume that $A \in L(H)$ and $X,Y \in L(K,H)$ with $R(A)$ closed and $R(X)\subseteq R(A),R(Y)\subseteq R(A^*)$.
We prove that
\begin{equation*}
(A-XY^*)^{+}=(I -(A^+XY^*)(A^+XY^*)^+)A^{+}(I-(XY^*A^+)^+(XY^*A^+))
\end{equation*}
if $XY^*A^+XY^*=XY^*$ and
\begin{equation*}
(A-XY^*)^{+}=(I -(A^+X)(A^+X)^+)A^{+}(I-(Y^*A^+)^+(Y^*A^+))
\end{equation*}
if $XY^*A^+X=X$ and $Y^*A^+XY^*=Y^*$.
These expressions generalize corresponding expressions of $(A-XY^*)^+$ given in  \cite{CHX} and
\cite{TFP}.

\section{preliminaries}
Let $T \in L(K,H)$, denote by $R(T)$ (resp. $N(T)$) the range (resp. kernal) of $T$.
Let $A \in L(H)$. Recall from \cite{WG} that $B \in L(H)$ is the Moore--Penrose inverse of $A$, if $B$ satisfies the following equations:
\begin{equation*}
ABA=A,~~BAB=B,~~(AB)^*=AB,~~(BA)^*=BA
\end{equation*}
In this case $B$ is denote by $A^+$. It is well--known $A$ has the Moore--Penrose inverse iff $R(A)$ is closed in $H$.
When $A^+$ exists, $R(A^+)=R(A^*),~N(A^+)=N(A^*)$ and $(A^+)^*=(A^*)^+$.

\begin{lemma}\label{YL}
Let $A \in L(H)$ with $R(A)$  closed and  $X,Y \in L(K,H)$\\
(1)~~$R(X)\subseteq R(A)$ iff $AA^+X=X$,~~$R(Y)\subseteq R(A^*)$ iff $Y^*A^+A=Y^*$.\\
(2)~~Suppose that $R(X)\subseteq R(A)$ and $R(Y)\subseteq R(A^*)$ then
\begin{equation*}
(A-XY^*)A^+(A-XY^*)=(A-XY^*)
\end{equation*}
iff $XY^*A^+XY^*=XY^*$.
\end{lemma}
\pf (1) Since $R(A)=R(AA^+)$ and $R(A^*)=R(A^+A)$, the assertion follows.

(2) Using (1), we can check directly that $(A-XY^*)A^+(A-XY^*)=(A-XY^*)$ if and only if
$XY^*A^+XY^*=XY^*$.

In order to compute $(A-XY^*)^+$, we need the following two lemmas which come from \cite{WXQ}.
\begin{lemma}\label{YLb}
Let $S \in L(H)$ be an idempotent operator. Denote by $O(S)$ the orthogonal projection of $H$ onto $R(S)$.
Then $I-S-S^*$ is invertible in L(H) and $O(S)=-S(I-S-S^*)^{-1}$.
\end{lemma}
\begin{lemma}\label{YLc}
Let $T,~B \in L(H)$ with $TBT=T$, Then $T^+=(I-O(I-BT))BO(TB)$.
\end{lemma}
\begin{lemma}\label{YLd}
Let $S\in L(H)$ be an idempotent operator. Then $O(S)=SS^+$ and $O(I-S)=I-S^+S$.
\end{lemma}
\pf $S^2=S$ implies that $R(S)$ is closed and $R(I-S)=N(S)=R(S^*)^\perp$. Thus, $S^+$ exists and
$O(S)=SS^+,\ O(I-S)=I-S^+S$.

\section{Main results}
\setcounter{equation}{0}

In this section, we will generalize Eq(\ref{EQ}) and Eq(\ref{EQ2}). Firstly, we have
\begin{proposition}\label{PP}
Let $A \in L(H)$ with $R(A)$ closed and  $X,Y \in L(K,H)$ with $R(X)\subseteq R(A)$ and $R(Y)\subseteq R(A^*)$.
Assume that $I-Y^*A^+X$ is invertible in L(H). Then $(A-XY^*)^+$ exists and
\begin{equation}
(A-XY^*)^+=A^++A^+X(I-Y^*A^+X)^{-1}Y^*A^+.
\end{equation}
\end{proposition}
\pf Put $B=A^++A^+X(I-Y^*A^+X)^{-1}Y^*A^+$. Simple computation shows that
$(A-XY^*)B=AA^+$ and $B(A-XY^*)=A^+A$ by Lemma \ref{YL} (1). Thus,
\begin{align*}
(A-XY^*)B(A-XY^*)&=A-XY^*,\quad\qquad B(A-XY^*)B=B,\\
((A-XY^*)B)^*&=(A-XY^*)B,\quad (B(A-XY^*))^*=B(A-XY^*),
\end{align*}
that is, $(A-XY^*)^+=B$.

Now we consider the case that $I-Y^*A^+X$ is not invertible, we have

\begin{theorem} \label{DL}
Let $A \in L(H)$ with $R(A)$ closed and  $X,Y \in L(K,H)$ with $R(X)\subseteq R(A)$ and
$R(Y)\subseteq R(A^*)$.
\begin{enumerate}
\item[(1)] If $XY^*A^+XY^*=XY^*$, then $(A-XY^*)^+$ exists and
\begin{equation}\label{4eq}
(A-XY^*)^{+}=(I -(A^+XY^*)(A^+XY^*)^+)A^{+}(I-(XY^*A^+)^+(XY^*A^+));
\end{equation}
Especially, if $XY^*A^+X=X$ and $Y^*A^+XY^*=Y^*$, then
\begin{equation}\label{5eq}
(A-XY^*)^{+}=(I -(A^+X)(A^+X)^+)A^{+}(I-(Y^*A^+)^+(Y^*A^+));
\end{equation}
\item[(2)] Assume that $R(A-XY^*)$, $R(A^+XY^*)$ and $R(XY^*A^+)$ are closed in $H$.
Then Eq(\ref{4eq}) implies that $XY^*A^+XY^*=XY^*$;
\item[(3)] Assume that $R(A-XY^*)$, $R(A^+X)$ and $R(Y^*A^+)$ are closed. Then
Eq(\ref{5eq}) indicates that $XY^*A^+X=X$ and $Y^*A^+XY^*=Y^*$.
\end{enumerate}
\end{theorem}
\pf (1) In this case, $(A-XY^*)A^+(A-XY^*)=(A-XY^*)$. Thus $R(A-XY^*)$ is closed,
i.e., $(A-XY^*)^+$ exists and hence
$$(A-XY^*)^{+}=(I-O(I-A^+(A-XY^*)))A^+O((A-XY^*)A^+)$$
by Lemma \ref{YL} (2). Since $(I-2A^+A)^2=I,~~(I-2A^+A)A^+=-A^+$,
\begin{align*}
A^+XY^*+(A^+XY^*)^*)&=(A^+XY^*+(A^+XY^*)^*)(2A^+A-I)\\
(I-A^+A)(I-A^+XY^*-(A^+XY^*)^*)&=I-A^+A.
\end{align*}
it follows that
\begin{align*}
O(I-A^+(A-XY^*))&=O(I-A^+A+A^+XY^*)\\
&=-(I-A^+A+A^+XY^*)(2A^+A-I-A^+XY^*-(A^+XY^*)^*)^{-1}\\
&=(I-A^+A+A^+XY^*)(I-2A^+A)(I-A^+XY^*-(A^+XY^*)^*)^{-1}\\
&=I-A^+A+O(A^+XY^*).
\end{align*}
Similarly, we also have
\begin{align*}
O((A-XY^*)A^+)&=-(A-XY^*)A^+(I-(AA^+-XY^*A^+)-(AA^+-XY^*A^+)^*)^{-1}\\
&=(-AA^++XY^*A^+)(I-2AA^++XY^*A^++(XY^*A^+)^*)^{-1}\\
&=(AA^+-XY^*A^+)(I-XY^*A^+-(XY^*A^+ )^*)^{-1}\\
&=AA^+-I+O(I-XY^*A^+).
\end{align*}
Therefore, we have
\begin{align*}
(A-XY^*)^{+}&=(I-O(I-A^+(A-XY^*)))A^+O((A-XY^*)A^+)\\
&=(A^+A-O(A^+XY^*))A^+O(I-XY^*A^+)\\
&=(I-O(A^+XY^*))A^+O(I-XY^*A^+).
\end{align*}

From $A^+XY^*A^+XY^*=XY^*$, we get that $A^+XY^*$ and $XY^*A^+$ are all idempotent operators.
It follow from Lemma \ref{YLd} that
$$
O(A^+XY^*)=(A^+XY^*)(A^+XY^*)^+,\quad O(I-XY^*A^+)=I-(XY^*A^+)^+(XY^*A^+).
$$
Therefore, we have
$$
(A-XY^*)^{+}=(I -(A^+XY^*)(A^+XY^*)^+)A^{+}(I-(XY^*A^+)^+(XY^*A^+)).
$$

When $XY^*A^+X=X$ and $Y^*A^+XY^*=Y^*$, we have $ R(A^+XY^*)=R(A^+X)$ and
$R(I-XY^*A^+)=N(Y^*A^+)$ so that
$$
O(A^+XY^*)=(A^+X)(A^+X)^+,\quad O(I-XY^*A^+)=I-(Y^*A^+)^+(Y^*A^+).
$$
and consequently, we get (\ref{5eq}).

(2) In this case,
$$
R((XY^*A^+)^*)=R((XY^*A^+)^+)\subseteq N((A-XY^*)^{+})=N((A-XY^*)^*),
$$
that is, $[N(XY^*A^+)]^\perp \subseteq [R(A-XY^*)]^\perp$. So $R(A-XY^*)\subseteq N(XY^*A^+)$
and consequently, $XY^*A^+XY^*=XY^*$.

(3) When Eq(\ref{5eq}) holds,
\begin{align*}
R((Y^*A^+)^*)=R((Y^*A^+)^+)&\subseteq N((A-XY^*)^{+})=N((A-XY^*)^*)\\
R((A-XY^*)^*)=R((A-XY^*)^+)&\subseteq N((A^+X)^+)=N((A^+X)^*).
\end{align*}
Then $R(A-XY^*)\subseteq N(Y^*A^+)$ and $R(A^+X)\subseteq N(A-XY^*)$. So
$$
Y^*A^+XY^*=Y^*,\quad XY^*A^+X=X.
$$

Suppose $H=\mathbb C^m$ and $K=\mathbb C^n$. Let $A\in L(H)$ and $X,\,Y\in L(K,H)$.
Since
\begin{align*}
\rank(A,X)=\rank A&\Leftrightarrow R(X) \subseteq R(A)\\
\rank\begin{pmatrix}A\\Y^*\end{pmatrix}=\rank A &\Leftrightarrow R(Y)\subseteq R(A^*),
\end{align*}
we can express Theorem \ref{DL} (1) as follows.
\begin{corollary}\label{Ca}
Let $A$ be an $m \times m$ matrix and $X,\,Y$ be two $m \times n$ matrices. Suppose that
$\rank(A,X)=\rank A$ and $\rank\begin{pmatrix}A\\Y^*\end{pmatrix}=\rank A$. Then
$$
(A-XY^*)^{+}=(I -(A^+XY^*)(A^+XY^*)^+)A^{+}(I-(XY^*A^+)^+(XY^*A^+))
$$
if $XY^*A^+XY^*=XY^*$ and
$$
(A-XY^*)^{+}=(I -(A^+X)(A^+X)^+)A^{+}(I-(Y^*A^+)^+(Y^*A^+))
$$
when $XY^*A^+X=X$ and $Y^*A^+XY^*=Y^*$.
\end{corollary}
Before ending this note, we give an example as follows.
\begin{example}Put $A=
\begin{pmatrix}
1&1&1&1\\0&0&1&1\\0&0&1&1\\0&0&0&1\end{pmatrix},\ X=\begin{pmatrix}1&1&0\\1&0&0\\1&0&0\\1&0&0
\end{pmatrix},\ Y=\begin{pmatrix}0&0&0\\0&0&0\\1&0&0\\1&0&1\end{pmatrix}.
$
Then
\begin{align*}
A^+=
\begin{pmatrix}
\frac{1}{2}&-\frac{1}{4}&-\frac{1}{4}&0\\\frac{1}{2}&-\frac{1}{4}&-\frac{1}{4}&0\\0&\frac{1}{2}&\frac{1}{2}&-1\\0&0&0&1
\end{pmatrix},&\quad XY^*=\begin{pmatrix}0&0&1&1\\0&0&1&1\\0&0&1&1\\0&0&1&1\end{pmatrix}\\
(A^+XY^*)^+=\begin{pmatrix}0&0&0&0\\0&0&0&0\\0&0&0&\frac{1}{2}\\0&0&0&\frac{1}{2}\end{pmatrix},
&\quad (XY^*A^+)^+=
\begin{pmatrix}
0&0&0&0\\ \frac{1}{4}&\frac{1}{4}&\frac{1}{4}&\frac{1}{4}\\ \frac{1}{4}&\frac{1}{4}&\frac{1}{4}&\frac{1}{4}\\0&0&0&0
\end{pmatrix}.
\end{align*}
It is easy to verify that $R(X)\subseteq R(A),~R(Y)\subseteq R(A^*)$ and $XY^*A^+XY^*=XY^*$.
So by Corollary \ref{Ca},
$
(A-XY^*)^{+}=
\begin{pmatrix}
\frac{1}{2}&0&0&0\\\frac{1}{2}&0&0&0\\0&0&0&-1\\0&0&0&0\end{pmatrix}.
$
\end{example}


\begin{thebibliography}{99}
\bibitem{WG} A. Ben-Israel and T.N.E. Greville, \emph{Generalized inverse:Theory and applications},
Wiley, New York, 1974.

\bibitem{WXQ}G. Chen and Y. Xue,\emph{The expression of generalized inverse of the perturbed
operators under type I perturbation in Hilbert spaces}, Linear Algebra Appl., 285(1998), 1--6.

\bibitem{CHX}Y. Chen, X. Hu and Q. Xu, \emph{The Moore-Penrose inverse of $A-XY^*$}, Journal of
Shanghai Normal University, 38(2009), 15--19

\bibitem{HS}H.V. Hsnderson and S.R. Searl, \emph{On deriving the inverse of a sum of matrices},
Siam Review, 23(1)(1981), 53--60.

\bibitem{KR}S. Kurt and A. Riedel, \emph{A Shermen-Morrison-Woodbury identity for rank augmenting}, matrices with application to centering.
Siam J. Math. Anal., 12(1)(1991), 80--95.

\bibitem{TFP}T. Steerneman and F.P. Kleij, \emph{Properties of the matrix
$A-XY^*$}, Linear Algebra Appl., 410(2005), 70--86.

\bibitem{WWH} W.W. Hager, \emph{Updating the inverse of a matrix}, Siam Review, 31(1989), 221--239.





\end{thebibliography}
\end{document}